\documentclass[reqno]{amsart}
\usepackage{amsmath,amssymb,amsthm,color}
\usepackage[colorlinks]{hyperref}
\usepackage{mathrsfs}
\usepackage{tikz}
\usetikzlibrary{fit}
\usetikzlibrary{shapes}
\usetikzlibrary{arrows}
\usetikzlibrary{decorations.markings}
\usetikzlibrary{positioning}

\newcommand{\colimit}{\varinjlim}
\newcommand{\ssets}{\textnormal{\textbf{SSet}}}
\newcommand{\CC}{\mathscr{C}}
\newcommand{\PP}{\mathscr{P}}
\newcommand{\iso}{\cong}
\newcommand{\tensor}{\otimes}
\newcommand{\homotopic}{\simeq}
\newcommand{\presheaves}[1]{\widehat{#1}}
\newcommand{\bottom}{\perp}
\newcommand{\cat}{\textnormal{\textbf{Cat}}}
\newcommand{\sets}{\textnormal{\textbf{Set}}}
\newcommand{\graphs}{\textnormal{\textbf{Graph}}}

\tikzset{mylabel/.style={font=\footnotesize}}
\tikzset{mymidlabel/.style={font=\footnotesize,fill=white}}

\newtheorem{theorem}{Theorem}[section]
\newtheorem{lemma}[theorem]{Lemma}
\newtheorem{proposition}[theorem]{Proposition}
\newtheorem{corollary}[theorem]{Corollary}

\theoremstyle{definition}
\newtheorem{definition}[theorem]{Definition}
\newtheorem{example}[theorem]{Example}

\theoremstyle{remark}

\begin{document}

\title{model structures from a monad on presheaves}
\author{Michael A. Warren}
\address{School of Mathematics\\Institute for Advanced Study\\Einstein Drive\\Princeton, New
  Jersey\\ 08540 USA}

\date{\today}

\begin{abstract}
  In this note we describe conditions under which the algebras
  for a monad on a presheaf category equipped with some additional
  structure are fibrant objects in a model structure.  We also prove
  that when these conditions are satisfied the resulting model
  structure is, in a suitable sense, the smallest model structure for
  which the units of the monad give a fibrant replacement.
\end{abstract} 

\maketitle

\section{Homotopy theory of presheaves}\label{section:cisinski}

Throughout $\CC$ is a small category and $\presheaves{\CC}$ denotes
the category of presheaves on $\CC$.  We recall some useful
results about the homotopy theory of presheaves
due to Cisinski  \cite{Cisinski:THDT,Cisinski:PCMTH}. 

\subsection{Elementary homotopy data}

\begin{definition}
  A \emph{cylinder} for a presheaf $X$ on $\CC$ consists of a
  tuple $(X\tensor I,\partial^{0},\partial^{1},\sigma)$ as
  indicated in the following diagram
  \begin{align*}
    \begin{tikzpicture}[auto]
      \node (leftX) at (0,1.5) {$X$};
      \node (rightX) at (3.5,1.5) {$X$};
      \node (cyl) at (1.75,1.5) {$X\tensor I$};
      \node (middleX) at (1.75,0) {$X$};
      \draw[->] (leftX) to node[mylabel] {$\partial_{0}$} (cyl);
      \draw[->] (rightX) to node[mylabel,swap] {$\partial_{1}$} (cyl);
      \draw[->] (leftX) to node[mylabel,swap] {$1_{X}$} (middleX);
      \draw[->] (rightX) to node[mylabel] {$1_{X}$} (middleX);
      \draw[->] (cyl) to node[mymidlabel,auto=false] {$\sigma$} (middleX);
    \end{tikzpicture}
  \end{align*}
  such that $[\partial^{0},\partial^{1}]:X+X\to X\tensor I$ is a
  monomorphism.

  A \emph{morphism of cylinders} $X\tensor I\to Y\tensor I$
  consists of maps $\phi:X\to Y$ and $\psi:X\tensor
  I\to Y\tensor I$ such that
  \begin{align*}
    \psi\circ\partial^{e} & = \partial^{e}\circ\phi\\
    \phi\circ\sigma & = \sigma\circ\psi
  \end{align*}
  for $e=0,1$.  This data determines a category
  $\text{Cyl}(\CC)$ of cylinders in $\presheaves{\CC}$.  

  A \emph{functorial cylinder} is a section of the forgetful functor
  $\text{Cyl}(\CC)\to \CC$.  

  \emph{Elementary homotopy data} on $\presheaves{\CC}$ consists of a
  functorial cylinder $(-\tensor I)$ which commutes with small colimits
  and preserves monomorphisms and for which each square
  \begin{align*}
    \begin{tikzpicture}[auto]
      \node (K) at (0,1.25) {$K$};
      \node (L) at (2.5,1.25) {$L$};
      \node (Kcyl) at (0,0) {$K\tensor I$};
      \node (Lcyl) at (2.5,0) {$L\tensor I$};
      \draw[->] (K) to node[mylabel] {$j$} (L);
      \draw[->] (K) to node[mylabel,swap] {$\partial^{e}$} (Kcyl);
      \draw[->] (L) to node[mylabel] {$\partial^{e}$} (Lcyl);
      \draw[->] (Kcyl) to node[mylabel,swap] {$j\tensor I$} (Lcyl);
    \end{tikzpicture}
 \end{align*}
  is a pullback when $j:K\to L$ is a monomorphism, for $e=0,1$.
\end{definition}
\begin{example}
  Multiplication by the simplicial interval $\Delta[1]$ provides the
  category $\ssets$ of simplicial sets with elementary homotopy data.
  We consider later a different choice of elementary homotopy data on $\ssets$.
\end{example}
Given elementary homotopy data $(-\tensor I)$ on $\presheaves{\CC}$ there is
an associated notion of homotopy obtained by defining a $(-\tensor
I)$\emph{-homotopy} from $f:X\to Y$ to $g:X\to Y$ to consist
of a map $\theta:X\tensor I\to Y$ such that
$\theta\circ\partial^{0}=f$ and $\theta\circ\partial^{1}=g$.  We write
$\homotopic$ for the homotopy relation and $\homotopic_{*}$ for the
equivalence relation on hom-sets generated by $\homotopic$.  $[X,Y]$
denotes the quotient of the hom-set $\ssets(X,Y)$ by
$\homotopic_{*}$.  Finally, note that these relations are
congruences for composition.

\subsection{Anodyne extensions}

Let $T$ be any set of maps in $\presheaves{\CC}$ and define
\begin{align*}
  \mathbf{\Lambda}(T) & := \{K\tensor I\cup L\tensor\partial I\to L\tensor
  I\;|\;K\to L\in T\}.
\end{align*}
Then, given elementary homotopy data consisting of a cylinder
$(-\tensor I)$ on
$\presheaves{\CC}$ and a set $S$ of monomorphisms we define a collection of
maps $\Lambda_{I}(S)$ inductively as follows:
\begin{align*}
  \Lambda_{I}^{0}(S) & := S\cup\{K\tensor I\cup
  L\tensor\{e\}\to L\tensor I\;|\;K\to L\in M,\;\;e=0,1\}
\end{align*}
and then:
\begin{align*}
  \Lambda_{I}^{n+1}(S) & :=
  \mathbf{\Lambda}\bigl(\Lambda^{n}_{I}(S)\bigr)\\
  \Lambda_{I}(S) & := \bigcup_{n\geq 0}\Lambda^{n}_{I}(S),
\end{align*}
where $M$ is any set of monomorphisms such that the saturated class
generated by $M$ is the class of all monomorphisms (such a $M$ will
always exist in a presheaf topos).  In
this case we define the class $\frak{A}_{I}(S)$ of anodyne extensions
generated by $I$ and $S$ by
\begin{align*}
  \frak{A}_{I}(S) & := ~^{\pitchfork}\!\bigl(\Lambda_{I}(S)^{\pitchfork}\bigr).
\end{align*}
The class $\frak{A}_{I}(S)$ is, in fact, the smallest saturated class
of monomorphisms which contains
\begin{align*}
  S\cup\{K\tensor I\cup L\tensor\{e\}\to L\tensor I\;|\;K\to L\in M,\;\;e=0,1\}
\end{align*}
and is closed under the operation $\mathbf{\Lambda}(-)$.  We say that
an arrow $f:X\to Y$ is a \emph{na\"{i}ve fibration} if it has the
RLP with respect to $\frak{A}_{I}(S)$.  A map $f:X\to Y$ is then
said to be a \emph{weak equivalence} if the induced map
\begin{align*}
  \begin{tikzpicture}
    \node (Left) {$[Y,A]$};
    \node (Right) [right=of Left] {$[X,A]$};
    \draw[->] (Left) to node[mylabel,auto] {$[f,A]$} (Right);
  \end{tikzpicture}
\end{align*}
is bijective for every (na\"{i}vely) fibrant object $A$.  Note that,
given elementary homotopy data $(-\tensor I)$ on $\presheaves{\CC}$, it is possible to
describe the sets $[Y,A]$ as consisting of the set of
maps $Y\to A$ modulo the equivalence relation generated by
$I$-homotopy.   We note the following lemma used by Cisinski in the
proof of Theorem \ref{theorem:cisinski} below.
\begin{lemma}
  \label{prop:fib_equiv}
  If $A$ is (na\"{i}vely) fibrant, then the homotopy relation
  $\homotopic$ is an equivalence relation on $\presheaves{\CC}(X,A)$, for
  any presheaf $X$.
\end{lemma} 
We will sometimes make free use of Lemma \ref{prop:fib_equiv} later on without
mention.  We now arrive at one of Cisinski's main results about the
homotopy theory of presheaf categories.
\begin{theorem}[Cisinski]\label{theorem:cisinski}
  If $I$ is elementary homotopy data on $\presheaves{\CC}$ and $S$ is a
  set of monomorphisms, then there is a model structure on
  $\presheaves{\CC}$ with the cofibrations the monomorphisms, the weak
  equivalences as above and the fibrations those maps having
  the RLP with respect to maps which are simultaneously cofibrations
  and weak equivalences. 
\end{theorem}
When $S$ is a set of monomorphisms and $I$ is elementary homotopy data
on $\presheaves{\CC}$ the model structure from Theorem
\ref{theorem:cisinski} is referred to as the \emph{$(S,I)$-generated
  model structure on $\presheaves{\CC}$} or, where $I$ is understood, as the
\emph{$S$-generated model structure on $\presheaves{\CC}$}.  Observe that in
such a model structure the $I$-homotopy equivalences are always also
weak equivalences.

\subsection{Example: Kan complexes}

When $\presheaves{\CC}$ is the category $\ssets$ of simplicial sets,
$I$ is the simplicial interval $\Delta[1]$ and the
set $S$ from above is the emptyset, Theorem \ref{theorem:cisinski}
gives the classical model structure on $\ssets$ in which weak
equivalences are weak homotopy equivalences and the fibrant objects
are exactly the Kan complexes.

\subsection{Example: Quasi-categories}

We include here as an example the case of the quasi-category model
structure on $\ssets$ \cite{Joyal:NOQC}.  Recall that here $J^{\infty}$ is the
\emph{infinite dimensional sphere}: the nerve of the groupoid interval
consisting of two object $\bottom$ and $\top$ with one non-identity isomorphism
$u:\bottom\to \top$ with inverse $d:\top\to \bottom$.  In this
instance $J^{\infty}$ is the appropriate interval for which we apply
Theorem \ref{theorem:cisinski}.

We write $\tau_{1}:\ssets\to \cat$ for the left-adjoint to the nerve
functor.  I.e., the left-Kan extension of the inclusion
$\mathbf{\Delta}\to \cat$ along the Yoneda embedding.  We say that
$\tau_{1}(X)$ is the \emph{fundamental category of $X$} and write
$\tau_{0}(X)$ for the set of isomorphism classes of objects of
$\tau_{1}(X)$.  A map of $f:X\to Y$ is said to be a \emph{weak
  categorical equivalence} if, for every quasi-category $A$, the
induced map
\begin{align*}
   \begin{tikzpicture}
    \node (Left) {$\tau_{0}(Y,A)$};
    \node (Right) [right=of Left] {$\tau_{0}(X,A)$};
    \draw[->] (Left) to node[mylabel,auto] {$\tau(f,A)$} (Right);
    \end{tikzpicture}
\end{align*}
is bijective.  Using a lifting property of $J^{\infty}$ established by
Joyal \cite{Joyal:QCKC} it is straightforward to establish that, when
$A$ is a quasi-category, the set $\tau_{0}(X,A)$ has another
description as the coequalizer
\begin{align*}
  \begin{tikzpicture}
    \node (Left) {$\ssets(X\times J^{\infty},A)$};
    \node (Mid) [right=of Left] {$\ssets(X,A)$};
    \node (Right) [right=of Mid] {$\tau_{0}(X,A)$.};
    \draw[transform canvas={yshift=0.5ex},->] (Left) to node[mylabel,auto] {$X^{\bottom}$} (Mid);
    \draw[transform canvas={yshift=-0.5ex},->] (Left) to
    node[mylabel,below] {$X^{\top}$} (Mid);
    \draw[->] (Mid) to (Right);
  \end{tikzpicture}
\end{align*}
Thus, the definition of weak categorical equivalences given here would
correspond exactly with that given by Theorem \ref{theorem:cisinski}
if it can be established that, for $J^{\infty}$ and an appropriate
set $S$ of maps as above, the quasi-categories are precisely the
na\"{i}vely fibrant objects.  

It follows, that the saturated class
generated by the set $S$ of inner horns
$\Lambda^{k}[n]\to \Delta[n]$ contains all maps of the form
\begin{align*}
  (\partial\Delta[n]\times J^{\infty})\cup(\Delta[n]\times
  \{e\})\to \Delta[n]\times J^{\infty}.
\end{align*}
Therefore, to see that $\text{sat}(S)=\frak{A}_{J^{\infty}}(S)$, it
suffices to prove that $\text{sat}(S)$ is closed under
$\mathbf{\Lambda}(-)$, which follows from results in \cite{Joyal:NOQC}.
Thus, the model structure coming from Theorem
\ref{theorem:cisinski} in this case is precisely the quasi-category
model structure.

\section{A model structure for $T$-algebras}

In this section we study model structures on $\presheaves{\CC}$
derived from a monad $T$ on $\presheaves{\CC}$ possessing appropriate
structure. I.e., we here describe conditions on a monad $T$ on
$\presheaves{\CC}$ which yield a model structure on $\presheaves{\CC}$
in which all $T$-algebras are fibrant.

We henceforth assume given a fixed choice of elementary homotopy data  $(-\tensor
I)$ for $\presheaves{\CC}$.

\subsection{$T$-weak equivalences}\label{section:we}

We arrive now at a modification of the notion of weak equivalence from above.
\begin{definition}
  A map $f:X\to Y$ is a $T$-\emph{weak equivalence (relative to
    $(-\tensor I)$)} if, for any
  $T$-algebra $A$, the induced map
  \begin{align*}
    \begin{tikzpicture}
      \node (Left) {$[Y,A]$};
      \node (Right) [right=of Left] {$[X,A]$};
      \draw[->] (Left) to node[mylabel,auto] {$[f,A]$} (Right);
    \end{tikzpicture}
  \end{align*}
  is a bijection. 
\end{definition}
It immediately follows from this definition that the $T$-weak
equivalences satisfy the familiar \emph{three-for-two} property and
that the collection $\frak{W}_{T}$ of $T$-weak equivalences contains
all $I$-homotopy equivalences.

\subsection{Axioms on $T$}

In this first attempt to arrive at a model structure on $\presheaves{\CC}$
related to $T$ in the way described above we first assume that $T$
satisfies the following two axioms.
\begin{description}
\item[(M1)] $T$ preserves monomorphisms and the units
  $\eta_{X}:X\to T(X)$ are monomorphisms.
\item[(M2)] There exists a set
  of presheaves $(G^{i})_{i\in I}$ on $\CC$ such that, for
  every presheaf $X$, the unit $\eta_{X}:X\to TX$ is in the
  saturated class generated by the maps
  \begin{align*}
    \begin{tikzpicture}
      \node (Left) {$G^{i}$};
      \node (Right) [right=of Left] {$T(G^{i})$};
      \draw[->] (Left) to node[mylabel,auto] {$\eta_{G^{i}}$} (Right);
    \end{tikzpicture}
  \end{align*}
  for $i\in I$.
\end{description}
Let $T$ be a monad on $\presheaves{\CC}$ which satisfies these axioms 
and write $\frak{G}$ for the generating set of units
$\{G^{i}\to TG^{i}\}$.  Then, by Cisinski's Theorem
\ref{theorem:cisinski}, there is a model structure on $\presheaves{\CC}$ in
which the weak equivalences, fibrations and
cofibrations are as described in Section \ref{section:cisinski} where now the
set $S$ is $\frak{G}$.  We now investigate conditions under which the
$T$-weak equivalences coincide with the weak equivalences obtained
from the $\frak{G}$-generated model structure.  
\begin{lemma}\label{lemma:units_acyclic_cof}
  All units $\eta_{X}:X\to T(X)$ are acyclic cofibrations in the
  $\frak{G}$-generated model structure on $\presheaves{\CC}$.
  \begin{proof}
    By the generation axiom, all units are in the saturated class
    $\text{sat}(\frak{G})$ of
    monomorphisms generated by $\frak{G}$.  Moreover, the collection
    of anodyne extensions $\frak{A}_{I}(\frak{G})$ is saturated and
    contains $\frak{G}$.  Therefore, $\text{sat}(\frak{G})$ is
    contained in $\frak{A}_{I}(\frak{G})$. 
  \end{proof}
\end{lemma}
\begin{lemma}
  \label{lemma:T_we_we}
  All $T$-weak equivalences are weak equivalences.
  \begin{proof}
    Assume $f:X\to Y$ is a $T$-weak equivalence and let a fibrant
    object $A$ be given.
    Then, by Lemma \ref{lemma:units_acyclic_cof}, there exists a retraction
    $\alpha:TA\to A$ as indicated in the following diagram:
    \begin{align*}
      \begin{tikzpicture}[auto]
        \node (LeftA) {$A$};
        \node (RightA) [right=of LeftA] {$A$.};
        \node (Bottom) [below=of LeftA] {$TA$};
        \draw[->] (LeftA) to node[mylabel] {$1_{A}$} (RightA);
        \draw[->] (LeftA) to node[mylabel,swap] {$\eta_{A}$} (Bottom);
        \draw[->,dashed] (Bottom) to node[mylabel,swap] {$\alpha$} (RightA);
      \end{tikzpicture}
    \end{align*}
    So, given a map $g:X\to A$ there exists, because
    $f$ is $T$-weak equivalence, a $k:Y\to TA$ such that
    $k\circ f$ is in the same homotopy class as $\eta_{A}\circ
    g$.  Defining $\bar{g}$ to be $\alpha\circ k$, it follows
    that $\bar{g}\circ f$ is in the same homotopy class as $g$.
    Moreover, this construction yields the same homotopy class given
    any $g':X\to A$ which is homotopic to $g$.  Thus, $f$ is a
    weak equivalence.
  \end{proof}
\end{lemma}

\noindent When $T$ satisfies \textbf{(M2)} it makes sense to consider
whether it also satisfies the following, rather strong, axiom:
\begin{description}
\item[(M3)] If $A$ is a $T$-algebra, then $A\to 1$ has the RLP
  with respect to elements of $\frak{A}_{I}(\frak{G})$.
\end{description}
Clearly \textbf{(M3)} is equivalent to the statement that all
$T$-algebras are fibrant and so, as might be expected, it is
in general non-trivial to establish that
it is satisfied by a given monad.  Consequently, the
following proposition is certainly in accord with the principle of
``preservation of difficulty''.
\begin{theorem}
  \label{theorem:main_1}
  Let $T$ be a monad on $\presheaves{\CC}$ satisfying
  \textnormal{\textbf{(M1)}-\textbf{(M3)}}, then there exists a model structure on
  $\presheaves{\CC}$ for which the weak equivalences are the $T$-weak
  equivalences, the cofibrations are the monomorphisms, all
  $T$-algebras are fibrant and the units $\eta_{X}:X\to T(X)$
  present $T(X)$ as a fibrant replacement of $X$.  Moreover, this
  model structure is precisely the one obtained by applying Theorem
  \ref{theorem:cisinski} to the set $\frak{G}$ of generating units. 
  \begin{proof}
    It suffices, by the foregoing lemmata, to prove that all weak
    equivalences are also $T$-weak equivalences.  That this is so
    follows immediately from \textbf{(M3)}.
  \end{proof}
\end{theorem}
When a monad $T$ on $\presheaves{\CC}$ with elementary homotopy data
$-\tensor I$ satisfies the hypotheses of Theorem \ref{theorem:main_1} we will refer
to the resulting model structure as the \emph{$T$-minimal model
  structure on $\presheaves{\CC}$ (with respect to $-\tensor I$)}.  We
now mention some consequences which justify the definition of $T$-weak
equivalences and which explain this choice of terminology.

\subsection{Further consequences of the axioms}

Throughout this section we assume that $T$ satisfies axioms
\textbf{(M1)}-\textbf{(M3)} unless otherwise stated.
\begin{lemma}
  \label{lemma:homotopy_preserved}
  Given maps $f,g:X\to Y$, $f\homotopic g$ implies $Tf\homotopic Tg$. 
  \begin{proof}
    We first prove that $f\homotopic g$ implies $Tf\homotopic Tg$, as
    this clearly implies the case for $\homotopic$.
    Suppose $\theta:X\tensor I\to Y$ is a homotopy from $f$ to
    $g$.  Then there exists a map
    \begin{align*}
      \begin{tikzpicture}
        \node (Left) {$(X\tensor I)\cup(TX\tensor\partial I)$};
        \node (Right) [right=of Left] {$TY$};
        \draw[->] (Left) to node[mylabel,auto] {$\theta'$} (Right);
      \end{tikzpicture}
    \end{align*}
    induced by the commutativity of the following square:
    \begin{align*}
      \begin{tikzpicture}[auto]
        \node (UL) {$X\tensor\partial I$};
        \node (UR) [right=of UL] {$TX\tensor\partial I$};
        \node (BL) [below=of UL] {$X\tensor I$};
        \node (BR) [right=of BL,below=of UR] {$TY$.};
        \draw[->] (UL) to node[mylabel] {$\eta_{X}\tensor
          1_{\partial I}$} (UR);
        \draw[->] (UL) to (BL);
        \draw[->] (UR) to node[mylabel] {$[Tf,Tg]$} (BR);
        \draw[->] (BL) to node[mylabel,swap] {$\eta_{Y}\circ\theta$} (BR);
      \end{tikzpicture}
    \end{align*}
    Thus, because $TY$ is fibrant, there exists a lift $\bar{\theta}$
    as indicated in the following diagram
    \begin{align*}
      \begin{tikzpicture}[auto]
        \node (UL) {$(X\tensor I)\cup(TX\tensor\partial I)$};
        \node (UR) [right=of UL] {$TY$.};
        \node (B) [below=of UL] {$TX\tensor I$};
        \draw[->] (UL) to node[mylabel] {$\theta'$} (UR);
        \draw[->] (UL) to (B);
        \draw[->,dashed] (B) to node[mylabel,swap] {$\bar{\theta}$} (UR);
      \end{tikzpicture}
    \end{align*}
   This $\bar{\theta}$ is the required homotopy.
  \end{proof}
\end{lemma}
\begin{corollary}
  \label{cor:induced_homotopic}
  Given a $T$-algebra $A$ and maps $f,g:X\to A$, if $f\homotopic
  g$, then $f'\homotopic g'$ where
  $f'$ is the induced $T$-algebra homomorphism
  $TX\to A$ and similarly for $g'$.
\end{corollary}
One important consequence of Lemma \ref{lemma:homotopy_preserved} is
the following alternative characterization of the weak equivalences.
\begin{lemma}\label{lemma:alternative_we}
  A map $f:X\to Y$ is a weak equivalence if and only if
  \begin{align*}
    \begin{tikzpicture}
      \node (L) {$TX$};
      \node (R) [right=of L] {$TY$};
      \draw[->] (L) to node[mylabel,auto] {$T(f)$} (R);
    \end{tikzpicture}
  \end{align*}
  is a homotopy equivalence with its homotopy inverse
  $\bar{f}:TY\to TX$ a $T$-algebra homomorphism.
  \begin{proof}
    Because, as remarked earlier, homotopy equivalences are always
    weak equivalences it suffices to prove the ``only if'' direction.

    First, assume $f$ is a weak equivalence.  Then
    $[f,TX]:[Y,TX]\to [X,TX]$ is a bijection and there exists a
    homotopy class 
    \begin{align*}
      \begin{tikzpicture}
        \node (L) {$[\;Y$};
        \node (R) [right=of L] {$TX\;]$};
        \draw[->] (L) to node[mylabel,auto] {$h$} (R);
      \end{tikzpicture}
    \end{align*}
    of maps $Y\to TX$ such that $[h\circ f]=[\eta_{X}]$.  Choose a
    representative $h:Y\to TX$ and then note that, because $TX$ is a
    $T$-algebra there exists a unique $T$-algebra homomorphism $\bar{f}:TY\to TX$
    such that
    \begin{align*}
      \begin{tikzpicture}[auto]
        \node (L) at (0,1.25) {$TY$};
        \node (M) at (1.25,0) {$Y$};
        \node (R) at (2.5,1.25) {$TX$};
        \draw[->,dashed] (L) to node[mylabel] {$\bar{f}$} (R);
        \draw[->] (M) to node[mylabel] {$\eta_{Y}$} (L);
        \draw[->] (M) to node[mylabel,swap] {$h$} (R);
      \end{tikzpicture}
   \end{align*}
    commutes.  By choice of $h$ we have that $hf\homotopic \eta_{X}$
    and therefore, by Corollary \ref{cor:induced_homotopic},
    \begin{align}
      \label{eq:use_cor}
      \bar{f}\circ T(f)\homotopic 1_{X}.
    \end{align}
    Also,
    \begin{align*}
      T(f)\circ \bar{f}\circ \eta_{Y}\circ f & = T(f)\circ\bar{f}\circ
      T(f)\circ\eta_{X}\\
      & \homotopic T(f)\circ\eta_{X}\\
      & = \eta_{Y}f,
    \end{align*}
    where the homotopy is by (\ref{eq:use_cor}).  But then
    \begin{align*}
      [f,TX]\bigl([T(f)\circ\bar{f}\circ\eta_{Y}\bigr) & = [f,TX]\bigl([\eta_{Y}]\bigr)
    \end{align*}
    and since $[f,TX]$ is injective
    $T(f)\circ\bar{f}\circ\eta_{Y}\homotopic\eta_{Y}$.  Therefore,
    applying Corollary \ref{cor:induced_homotopic}, 
    \begin{align*}
      T(f)\circ\bar{f} & \homotopic 1_{Y}.
    \end{align*}
  \end{proof}
\end{lemma}
The following proposition together with its corollary show that the
model structure on $\presheaves{\CC}$ established in Proposition
\ref{theorem:main_1} above is, in a suitable sense, the least model
structure on $\presheaves{\CC}$ for which the units $\eta_{X}:X\to TX$
are fibrant replacements.
\begin{proposition}\label{prop:justifies_we}
  If $\frak{W}'$ is any collection of maps in $\presheaves{\CC}$ 
  such that
  \begin{enumerate}
  \item $\frak{W}'$ has the 3-for-2 property;
  \item $\frak{W}'$ contains all homotopy equivalences; and
  \item $\frak{W}'$ contains all units $\eta_{X}:X\to T(X)$,
  \end{enumerate}
  then the collection $\frak{W}$ of weak equivalences is contained in $\frak{W}'$.  
  \begin{proof}
    Let a map $f:X\to Y$ in $\frak{W}$ be given.  Then, by Lemma
    \ref{lemma:alternative_we}, there exists a
    homotopy inverse $\bar{f}:TY\to TX$ of $Tf$ which is also a
    $T$-algebra homomorphism.  Then $Tf$ is in $\frak{W}'$ by
    assumption and, because
    \begin{align*}
      \begin{tikzpicture}[auto]
        \node (UL) {$X$};
        \node (UR) [right=of UL] {$Y$};
        \node (BL) [below=of UL] {$TX$};
        \node (BR) [right=of BL,below=of UR] {$TY$};
        \draw[->] (UL) to node[mylabel] {$f$} (UR);
        \draw[->] (UL) to node[mylabel,swap] {$\eta_{X}$} (BL);
        \draw[->] (BL) to node[mylabel,swap] {$Tf$} (BR);
        \draw[->] (UR) to node[mylabel] {$\eta_{Y}$} (BR);
      \end{tikzpicture}
    \end{align*}
    commutes, the 3-for-2 axiom for $\frak{W}'$ implies that $f$ is
    also in $\frak{W}'$.
  \end{proof}
\end{proposition}

\section{Examples}

We now give two simple examples to illustrate that there do in fact
exist monads and elementary homotopy data for which the
conditions \textbf{(M1)} through \textbf{(M3)} are satisfied.

\subsection{Monoids}

Consider the free-monoid monad $T$ on the category $\sets$ of sets.
We will prove that this monad satisfies conditions \textbf{(M1)} to
\textbf{(M3)}.  First, observe that \textbf{(M1)} is trivially
satisfied.  Next, we take as our set of generating presheaves the set
$\mathbb{N}$ of natural numbers (coded as, say, Zermelo ordinals).  We
will now prove that the units of the natural numbers give rise to all
units.  Let $\frak{G}$ be the set of units of natural numbers.
\begin{lemma}
  \label{lemma:monoids_M2}
  With this $\frak{G}$, \textnormal{\textbf{(M2)}} is
  satisfied.
  \begin{proof}
    Let a set $X$ be given and denote by $\PP_{f}(X)$ the set
    of all pairs $(S,\sigma)$ such that $S$ is a finite subset of $X$
    and $\sigma$ is an isomorphism $|S|\iso S$ of sets, where $|S|$
    denotes the cardinal number of $S$.  The map 
    \begin{align*}
      \begin{tikzpicture}
        \node (L) at (0,0) {$\coprod_{(S,\sigma)\in\PP_{f}(X)}|S|$};
        \node (R) at (4.5,0)
        {$\coprod_{(S,\sigma)\in\PP_{f}(X)}T|S|$};
        \draw[->] (L) to node[mylabel,auto] {$\coprod_{(S,\sigma)}\eta_{|S|}$} (R);
      \end{tikzpicture}
    \end{align*}
    is in the saturated class generated by $\frak{G}$, where $|S|$
    denotes the cardinality of $S$.  We claim that there exist
    maps making
    \begin{align*}
      \begin{tikzpicture}[auto]
        \node (UL) at (0,2) {$X$};
        \node (UM) at (2,2) {$\coprod_{(S,\sigma)}|S|$};
        \node (UR) at (4,2) {$X$};
        \node (BL) at (0,0) {$TX$};
        \node (BM) at (2,0) {$\coprod_{(S,\sigma)}T|S|$};
        \node (BR) at (4,0) {$TX$};
        \draw[->] (UL) to node[mylabel,swap] {$\eta_{X}$} (BL);
        \draw[->] (UR) to node[mylabel] {$\eta_{X}$} (BR);
        \draw[->] (UM) to node[auto=false,mymidlabel] {$\coprod_{(S,\sigma)}\eta_{|S|}$}
        (BM);
        \draw[->] (UL) to node[mylabel] {$s$} (UM);
        \draw[->] (UM) to node[mylabel] {$r$} (UR);
        \draw[->] (BL) to node[mylabel,swap] {$u$} (BM);
        \draw[->] (BM) to node[mylabel,swap] {$v$} (BR);
      \end{tikzpicture}
    \end{align*}
    a retract diagram.  Given $x\in X$, let $s(x):=((\{x\},!),0)$ and let
    $r$ be the canonical map determined by the maps
    \begin{align*}
      \begin{tikzpicture}
        \node (L) {$|S|$};
        \node (M) [right=of L] {$S$};
        \node (R) [right=of M] {$X$};
        \draw[->] (L) to node[mylabel,auto] {$\sigma$} (M);
        \draw[->] (M) to (R);
      \end{tikzpicture}
    \end{align*}
    for $(S,\sigma)$ in $\PP_{f}(X)$, where the unnamed map is
    the inclusion.

    Every word $w=x_{1}\cdots x_{n}$ of length $n$ in $TX$ determines
    a finite subset $S=\{x_{1},\ldots,x_{n}\}$ of $X$ together the
    obvious isomorphism $\sigma:|S|\to S$ and we let
    $u(w):=((S,\sigma),(\sigma(x_{1})\cdots \sigma(x_{n})))$.
    Finally, $v$ is induced by the map $T|S|\to TX$, for
    $(S,\sigma)$ in $\PP_{f}(X)$, obtained by sending a word
    $(a_{1}\cdots a_{n})$ to the word
    $(\sigma(a_{1})\cdots\sigma(a_{n}))$.  This data clearly
    determines a retract diagram as claimed and therefore $\eta_{X}$
    is in the saturated class generated by $\frak{G}$.
  \end{proof}
\end{lemma}
Next, note that $\sets$ possesses a natural choice of
elementary homotopy data given by cartesian product with the subobject
classifier $2=\{0,1\}$ (using the subobject classifier in this way is
typical and arises in \cite{Cisinski:THDT} as well).  This is the
elementary homotopy data which we will employ.  With this choice we
have
\begin{lemma}
  \label{lemma:monoids_M3}
  With these choices, \textnormal{\textbf{(M3)}} is satisfied.
  \begin{proof}
    Let a monoid $M$ be given.  Since $M$ has the RLP with respect to
    all elements of $\frak{G}$, it suffices to show that there exists
    a map $\bar{f}$ making 
    \begin{align*}
      \begin{tikzpicture}[auto]
        \node (UL) {$K\times 2\cup L\times\{e\}$};
        \node (UR) [right=of UL] {$M$};
        \node (BL) [below=of UL] {$L\times 2$};
        \draw[->] (UL) to node[mylabel] {$f$} (UR);
        \draw[->] (UL) to (BL);
        \draw[->,dashed] (BL) to node[mylabel,swap] {$\bar{f}$} (UR);
      \end{tikzpicture}
    \end{align*}
    commute, for any such $f$ ($e=0,1$) and any monomorphism $K\to L$.
    The existence of $\bar{f}$ is trivial in this case (and does not
    even require $M$ to be a monoid).
  \end{proof}
\end{lemma}

\subsection{Categories}

Let $\cat$ denote the category of small categories.  Let $\mathbf{L}$
denote the free category on the graph consisting of two distinct objects $0$ and $1$ and
two distinct edges $u:0\to 1$ and $d:1\to 0$.

Let $\graphs$ be the category of directed graphs and let $T$ be the
free category monad on $\graphs$.  As in the case of monoids, $T$
trivially satisfies condition \textbf{(M1)}.  In what follows, we
denote by $[n]$ the finite ordinal $\{0,1,\ldots,n\}$ regarded as a
category.  We will prove that when $\frak{G}$ is the set of units of
finite ordinals $[n]$ for $n\geq 0$, the generation axiom is
satisfied.
\begin{lemma}
  \label{lemma:cats_M2}
  The free category monad $T$ on $\graphs$ satisfies condition
  \textnormal{\textbf{(M2)}}.
  \begin{proof}
    Let $\frak{G}$ be as just described and assume that $G$ is an
    arbitrary graph.  We first form a new graph $\tilde{G}(0)$ as
    the pushout
    \begin{align*}
      \begin{tikzpicture}[auto]
        \node (UL) at (0,1.5) {$\coprod_{[0]\to G}[0]$};
        \node (UR) at (2.5,1.5) {$G$};
        \node (BL) at (0,0) {$\coprod_{[0]\to G}T[0]$};
        \node (BR) at (2.5,0) {$\tilde{G}(0)$};
        \draw[->] (UL) to (UR);
        \draw[->] (UL) to node[mylabel,swap] {$\coprod\eta_{[0]}$} (BL);
        \draw[->] (BL) to (BR);
        \draw[->] (UR) to node[mylabel] {$h_{0}$} (BR);
      \end{tikzpicture}
    \end{align*}
    and we observe that there is a canonical induced map
    $k_{0}:\tilde{G}(0)\to TG$ such that $k_{0}\circ h_{0}=\eta_{G}$
    and $\coprod T[0]\to\tilde{G}(0)\to TG$ is the induced map
    $\coprod T[0]\to TG$.  In general, assuming we have constructed
    $\tilde{G}(n)$ we form the pushout
    \begin{align*}
      \begin{tikzpicture}[auto]
        \node (UL) at (0,1.5) {$\coprod_{[n+1]\to\tilde{G}(n)}[n+1]$};
        \node (UR) at (3.5,1.5) {$\tilde{G}(n)$};
        \node (BL) at (0,0) {$\coprod_{[n+1]\to\tilde{G}(n)}T[n+1]$};
        \node (BR) at (3.5,0) {$\tilde{G}(n+1)$};
        \draw[->] (UL) to (UR);
        \draw[->] (UL) to node[mylabel,swap] {$\coprod\eta_{[n+1]}$} (BL);
        \draw[->] (BL) to (BR);
        \draw[->] (UR) to node[mylabel] {$h_{n+1}$} (BR);
      \end{tikzpicture}
    \end{align*}
    and we obtain $k_{n+1}:\tilde{G}(n+1)\to TG$ with $h_{n+1}\circ
    k_{n+1}=k_{n}$ and also commuting with the map $\coprod
    T[n+1]\to\tilde{G}(n+1)$ in the required way.  This gives us a
    countable sequence $G\to\tilde{G}(0)\to\tilde{G}(1)\to\ldots$ and
    we take the colimit $\tilde{G}=\colimit_{n}\tilde{G}(n)$.  Since
    each of the maps $h_{n}$ is in the saturated class generated by
    $\frak{G}$ it follows that the resulting map $h:G\to\tilde{G}$ is
    also in this class.  There is also a canonical map $k:\tilde{G}\to
    TG$ such that $k\circ h=\eta_{G}$. 

    Since the saturated class generated by $\frak{G}$ is closed under
    retracts it suffices to define a section $s:TG\to \tilde{G}$ of
    the map $k$.  On vertices $s$ is just the identity.  An edge $a\to
    b$ in $TG$ is a sequence $(f_{1},f_{2},\ldots,f_{n})$
    of composable edges from $G$.  (Note that such a sequence is
    allowed to be empty $()$ when $a=b$.)  This gives graph
    homomorphisms $[n]\to G$ and, composing with the maps
    $h_{0},h_{1},\ldots,h_{n-1}$, $[n]\to\tilde{G}(n-1)$.  Such a sequence possesses a
    composite $(f_{n}\circ\cdots\circ f_{1})$ in $\tilde{G}(n)$
    and therefore determines an edge $[f_{n}\circ\cdots\circ
    f_{1}]:a\to b$ in $\tilde{G}$ and we take $s(f_{1},\ldots,f_{n})$
    to be this edge (here the square brackets indicate
    that this edge is actually itself an equivalence class).  With
    this definition, $s$ is clearly a section of $k$.
  \end{proof}
\end{lemma}
In $\graphs$ we have elementary homotopy data given by $(-\times
I)$ where $I$ is the graph with two distinct objects $0$ and $1$ and
two distinct edges $u:0\to 1$ and $d:1\to 0$.
\begin{lemma}
  \label{lemma:cats_M3}
  Condition \textnormal{\textbf{(M3)}} is satisfied for this choice of
  elementary homotopy data.
  \begin{proof}
    Given 
    \begin{align*}
      \begin{tikzpicture}[auto]
        \node (UL) {$K\times I\cup L\times\{0\}$};
        \node (UR) [right=of UL] {$C$};
        \node (BL) [below=of UL] {$L\times I$};
        \draw[->] (UL) to node[mylabel] {$f$} (UR);
        \draw[->] (UL) to (BL);
      \end{tikzpicture}
    \end{align*}
    with $C$ a category, we define $\bar{f}:L\times I\to C$ on
    vertices by setting
    \begin{align*}
      \bar{f}(x,t) & :=
      \begin{cases}
        f(x,0) & \text{ when }x\notin K\\
        f(x,t) & \text{ otherwise.}
      \end{cases}
    \end{align*}
    For edges, we will consider each of the different cases one at a
    time.  Assume given an edge $(e,j):(a,s)\to(b,t)$ in $L\times I$.
    For such $(e,j)$ in the domain of $f$ we simply define
    $\bar{f}(e,j):=f(e,j)$.  As such, we now focus on those cases
    where $(e,j)$ is not in the domain $\text{dom}(f)$ of $f$.
    \begin{description}
    \item[$a\notin\text{dom}(f)$ and $b\notin\text{dom}(f)$] If $e=a=b$, then
      $\bar{f}(e,j):=1_{f(e,0)}$.  Otherwise, $\bar{f}(e,j):=f(e,0)$.
    \item[$a\in \text{dom}(f)$ and $b\notin \text{dom}(f)$] If $j=1$, then
      $\bar{f}(e,j):=f(e,0)\circ f(a,d)$.  If $j=d$, then
      $\bar{f}(e,j):=f(e,0)\circ f(a,d)$.  If $j=u$, then $\bar{f}(e,j):=f(e,0)$.
    \item[$a\notin \text{dom}(f)$ and $b\in \text{dom}(f)$]
      Dual to the previous case.
    \item[$a\in \text{dom}(f)$ and $b\in \text{dom}(f)$]
      If $j=1$, then $a,b\in K$ and $\bar{f}(e,j):=f(b,u)\circ
      f(e,0)\circ f(a,d)$.  If $j=d$, then $\bar{f}(e,j):=f(e,0)\circ
      f(a,d)$.  If $j=u$, then $\bar{f}(e,j):=f(b,u)\circ f(e,0)$.
   \end{description}
   Thus, we have defined $\bar{f}$.  Similarly, a category $C$
   has the RLP with respect to those maps $K\times I\cup
   L\times\{1\}\to L\times I$ for $K\to L$ a monomorphism.
  \end{proof}
\end{lemma}

\subsection*{Acknowledgements}

Although these results were originally obtained a number of years ago,
the final preparation of this note was made while the author
received support from the Oswald Veblen Fund, for which he is grateful.

\newcommand{\SortNoop}[1]{}


\begin{thebibliography}{10}

\bibitem{Cisinski:THDT}
Denis-Charles Cisinski, Th{\'{e}}ories homotopiques dans les topos,
  \emph{J. Pure Appl. Algebra} \textbf{174}, 2002, 43--82.

\bibitem{Cisinski:PCMTH}
\bysame, \emph{Les pr{\'{e}}faisceaux comme mod{\`{e}}les des types
  d'homotopie}, Ast{\'{e}}risque, vol.~308, Soc. Math. France, 2006.

\bibitem{Joyal:QCKC}
Andr{\'{e}} Joyal, Quasi-categories and {K}an complexes, \emph{J. Pure
  Appl. Algebra} \textbf{175}, 2002, 207--222.

\bibitem{Joyal:NOQC}
\bysame, The theory of quasi-categories and its applications.
\newblock In {\em Advanced Course on Simplicial Methods in Higher
  Categories}, vol.~2, pages 149--496, Centre de Recerca
Matem\`{a}tica, Barcelona, 2008.

\end{thebibliography}
\end{document}